\documentclass[12pt,onecolumn]{elsart3p}

\usepackage{amssymb}
\usepackage{bm}
\usepackage{graphicx}
\usepackage{amsfonts}
\usepackage{color}
\usepackage{epstopdf}
\usepackage[square,comma,sort,numbers]{natbib}
\usepackage{multirow}
\usepackage{subfigure}
\usepackage{layout}

\def\fg#1{Fig. \ref{#1}} 
\def\eq#1{Eq. (\ref{#1})} 

\def\eqn#1{\begin{equation} \label{#1}}
\def\enn{\end{equation}}
\def\eqN{\[}
\def\enN{\]}

\def\e#1{eq (\ref{#1})}

\def\s#1{Section \ref{#1}}

\begin{document}
 \bibliographystyle{plainnat}

\begin{frontmatter}

\title{On the acceleration of spatially distributed agent-based computations: \\
a patch dynamics scheme}

\author[a,b]{Ping Liu},
\author[c]{Giovanni Samaey},
\author[d]{C. William Gear},
\author[d,e]{Ioannis G. Kevrekidis,\corauthref{abc} },
\address[a]{Department of Molecular, Cellular and Developmental Biology, Yale University, West Haven, CT 06516, USA.}
\address[b]{Systems Biology Institute, Yale University, West Haven, CT 06516, USA.}
\address[c]{Department of Computer Science, KU Leuven, 3001 Leuven, Belgium}
\address[d]{Department of Chemical and Biological Engineering, Princeton University, Princeton, NJ 08544, USA.}
\address[e]{Program in Applied and Computational Mathematics, Princeton University, Princeton, NJ 08544, USA.}
\corauth[abc]{corresponding author: Ioannis G. Kevrekidis, Email address: yannis@arnold.princeton.edu, Phone: +1-609 258-2818. Fax: +1-609 258-0211.}

\begin{abstract}

In recent years, individual-based/agent-based modeling has been applied
to study a wide range of applications, ranging from engineering problems to phenomena in sociology, economics and biology.
Simulating such agent-based models over extended spatiotemporal domains can be
prohibitively expensive due to stochasticity and the presence of multiple scales.
Nevertheless, many agent-based problems exhibit
smooth behavior in space and time on a macroscopic scale, suggesting that a
useful coarse-grained continuum model could be obtained.
%
For such problems, the equation-free framework \cite{G:01,G:02,G:03} 
can significantly reduce the computational cost.
Patch dynamics is an essential component of this framework. 
This scheme is designed to perform numerical simulations of an unavailable
macroscopic equation on macroscopic time and length scales; it uses appropriately
initialized simulations of the fine-scale agent-based model in a number
of small ``patches'', which cover only a fraction of the spatiotemporal domain.
%
In this work, we construct a finite-volume-inspired
\emph{conservative} patch dynamics scheme and apply it to a
financial market agent-based model based on the work of Omurtag and Sirovich
\cite{Omurtag:01}. 
%
We first apply our patch dynamics scheme to a continuum
approximation of the agent-based model, to study its performance and analyze its accuracy.
We then apply the scheme to the agent-based model itself.
Our computational experiments indicate that here, typically,
the patch dynamics-based simulation requires only $20\%$ of the full
agent-based simulation in space, and need occur over only
$10\%$ of the temporal domain.
\end{abstract}

\end{frontmatter}

\section{Introduction}

Over the last years agent-based modeling has become a powerful simulation
technique for a number of applications, varying from ecology
\cite{Ecology:01,Ecology:02,Ecology:03} to economics and financial systems
\cite{Finance01,Finance02,Finance03}, and from traffic and supply chain networks
\cite{Traffic:01,Traffic:02,Traffic:03} to biological \cite{Bio:01,Bio:02,Bio:03}
and social \cite{Social:01,Social:02,Social:03} systems.
Agents are treated as unique and discrete entities, allowing a straightforward way to incorporate detailed interactions between
individuals via a set of microscopic rules.
Modelers use the fine-scale evolution at the individual agent level to
study the macroscopic, population level dynamical behavior, often referred to as \emph{emergent} behavior. It would of course be desirable (and computationally much more efficient) to simulate this emergent behavior through a continuum macroscopic model. However, accurately accounting for the consequences of individual interactions at the macroscopic level can be highly nontrivial. As a consequence, accurate macroscopic equations are typically unavailable, even in cases when a clear separation of time scales strongly suggests that such equations \emph{should} exist.
In this work, we will consider a financial market agent-based model, in which a stochastic
process describes the ``motion'' of an individual agent along a spatial domain representing the agent's propensity to buy or sell.
This stochastic process needs to be simulated for a large number of (interacting) agents over long times.
Doing this over the entire spatiotemporal domain of interest to observe the
relevant macroscopic variables (e.g., agent density) incurs a
prohibitive computational cost.

The recently developed \emph{equation-free framework} (\cite{G:01,G:02,G:03}) can be applied to exploit
the time-scale separation between the (fine-scale) individual dynamics and the (coarse-scale) population-level behavior by performing the expensive agent-based simulations in a grid of small patches, which cover only a fraction
of the spatiotemporal domain. Exploiting smoothness, the solution for the entire region is  approximated by repeated interpolation/extrapolation of the macroscopic variables recorded within these patches.

This framework is built around the central idea of a coarse time-stepper, which advances the coarse variables over a time interval of size $\delta t$. It consists of the following steps: (1) \emph{lifting}, i.e., the creation
of appropriate initial conditions for the microscopic model; (2) fine-scale \emph{evolution}; and (3) \emph{restriction},
i.e., the estimation of macroscopic observables from the fine-scale solution. This coarse time-stepper can subsequently be used as input for
time-stepper-based ``wrapper'' algorithms performing macroscopic numerical analysis tasks.
These include, for example, time-stepper based bifurcation codes to perform
coarse  bifurcation analysis for the unavailable macroscopic equation \cite{bifur1,bifur2,bifur3},
and other system level tasks such as rare event analysis \cite{rare1},
control \cite{control1,control2} and optimization \cite{optim1}.

The patch dynamics scheme \cite{G:03} is a combination of
two sub-schemes - the gap-tooth scheme (interpolate macroscopic properties in
space) \cite{Giovanni2005} and the coarse projective integration scheme (extrapolate
macroscopic properties in time) \cite{Proj2004}. First, in the gap-tooth scheme
\cite{Giovanni2005}, microscopic simulations are only performed in a number of
small subdomains of the spatial domain (intervals for the case of one space dimension); in-between these patches lie what
we call gaps.
We continue the discussion for the one-space-dimension paradigm that fits our application.
A coarse time-$\delta t$ map is then constructed as follows. We first
choose a number of macroscopic grid points and a small interval (``tooth")) around
each grid point; we initialize the fine scale consistently with a macroscopic initial condition profile; apply the microscopic solver within each
interval, using appropriate boundary conditions for time $\delta t$; and, finally, obtain macroscopic
information at time $\delta t$ (e.g., by computing the average density in each tooth).
In \cite{Giovanni2005}, the initial conditions for each of these teeth are obtained via interpolation over neighboring teeth; this can be seen to be equivalent to constructing a macroscopic finite difference approximation. In the gap-tooth scheme, this procedure is then repeated. We refer to \cite{gaptooth2003} for an illustration of this scheme with particle-based simulations of the viscous
Burgers equation.

Second, exploiting the smoothness of coarse variable evolution in the {\em temporal} domain,
one can combine the gap-tooth scheme with a projective integration scheme
\cite{Proj2004}. In this context, we perform a number of gap-tooth steps of size $\delta t$ to
obtain an estimate of the \emph{time derivative} of the unavailable macroscopic equation.
Based on this estimate, a projective step is subsequently taken over time  $\Delta t \gg \delta t$. This combination has been termed {\em patch
dynamics} \cite{G:03}.

In previous work, the gap-tooth and patch dynamics
schemes  have been applied to study a diffusion homogenization
problem \cite{Giovanni2005,Giovanni2006}, as well as a biological dispersal model \cite{Radek-PhysD-2006}. For the kind of agent-based problems
we study in this work, a conservation law needs to be satisfied. For this reason,
instead of basing the gap-tooth scheme on a finite difference approximation, we design a finite volume-based patch
dynamics scheme. We first apply this scheme to an approximate continuum model of the agent-based dynamics to illustrate the scheme, as well as to facilitate our error analysis. In this case, our ``inner'' microscopic
model is a fine discretization of this continuum model, and our microscopic simulator is a classical finite volume scheme within each tooth. In general, a given
microscopic code only allows us to run it with a set of predefined boundary
conditions. It is highly non-trivial to impose macroscopically inspired
boundary conditions on such microscopic codes, see, e.g. \cite{BCctrl} for
a control-based strategy.  As in \cite{Giovanni2006}, the boundary conditions on the teeth are imposed {\em via the use of buffer regions}, surrounding the teeth,  to ``protect" the tooth inside each simulation unit from boundary
artifacts. We call ``simulation unit" the union of a tooth and its edge buffers.
After several gap-tooth steps, we project the macroscopic properties forward in
time for a bigger step $\Delta t$. Based on the projected macroscopic properties, we again construct consistent microscopic initial conditions for the simulation units and the cycle repeats. We
study the performance of this  patch dynamics scheme,
and compute its order of accuracy both analytically and numerically.
Subsequently, we apply the scheme to the agent-based
model, and observe that, typically, this requires
the expensive microscopic simulations in only $20\%$ of the spatial domain and $10\%$ of the temporal domain.

These micro/macro ideas also form the basis of the heterogeneous multiscale method \cite{HMM2003,HMM2007}. There, a macro-scale solver is combined with an estimator for quantities that are unknown because the macroscopic equation is not available.
This estimator subsequently uses appropriately constrained runs of the
microscopic model.  For a recent overview of results obtained in that framework, we refer to \cite{HMM2013-ActaNUmerica}.

The remainder of this paper is organized as follows: in section 2 we introduce an agent-based
financial market model based on the work of \cite{Omurtag:01} and the associated continuum approximation they derive.
In section 3 we
describe the ``inner'' finite volume scheme for the continuum model. In section 4 we
describe the patch dynamics scheme. In section 5 we first present patch
dynamics results using the discretization introduced in section 3 as the inner solver, and then analyze the order of
consistency of this scheme both analytically
and numerically. The main result, the agent-based patch dynamics computations, are also presented in section 5.
We conclude with a brief summary and discussion in section 6.

\section{The models}

\subsection{The agent-based model} \label{sec:agent}
 We consider a financial market model initially
 described by Omurtag and Sirovich \cite{Omurtag:01}. This model
 simulates the actions of buying and selling by a large population of $N$
 interacting individuals in the presence of mimesis. In this model, the
 $i$-th agent's propensity to buy or sell is indicated by its state $x_i(t) \in (-1,1)$
which
 evolves 
 according to two coupled processes. The first process is the exponential decay, at constant rate $\gamma$,
 of $x_i(t)$ towards the neutral state $0$.
 This decay implies that each agent gradually forgets
 its current state and tends to eventually become neutral in its preference for
 buying or selling.

The second process is a stochastic process denoted by $I_i(t)$ which represents
the effect of incoming information to each agent's state $x_i(t)$.
$I_i(t)$ incorporates three factors: (1) the arrival times of incoming information; (2) the type of this information (i.e., ``good'' or ``bad''); and (3) the quantitative effect of this information on $x_i(t)$. The arrival of incoming
information is assumed to follow a Poisson distribution with mean arrival
frequency $(\nu^+ + \nu^-)$, where $\nu^+$ denotes the mean arrival
frequency for good news and $\nu^-$ denotes the one for bad news. The mean
arrival frequencies are given by
\begin{equation}
                   \nu^{\pm} = \nu_{ex}^{\pm}+gR^{\pm},
\label{nu:eqn}
\end{equation}
where the parameters $\nu_{ex}^+$ and $\nu_{ex}^-$ represent the contribution
of \emph{external} information each individual receives from its environment (e.g., mass media news or opinions of stock market consultants) and are assumed time-independent. The quantities $R^+$ and $R^-$ are
the \emph{buy rates} and \emph{sell rates}, respectively, defined as the number of buys or sells
happening in the market per unit time over a small finite time horizon (which
we call \emph{reporting horizon}) normalized by the total number of agents in the
market. The parameter $g$ is a feedback constant which quantifies the extent
to which the individuals' propensity to buy or sell is determined by the buy and sell
rates. Note that the buy and sell rates are collective properties of the population
(i.e., all the agents in the market), which implies that the second term of (\ref{nu:eqn}),
$gR^{\pm}$, embodies how each individual agent's state is affected
by the collective behavior of the entire population at that moment in time.

Each arriving ``quantum'' of information  has probability $\nu^+/(\nu^+ + \nu^-)$ to be good
news and $\nu^-/(\nu^+ + \nu^-)$ to be bad news. When good news arrives for agent $i$, the value of $x_i(t)$ is increased
instantaneously - jumps - by a fixed amount $\epsilon^+>0$; similarly,
when bad news arrives, $x_i(t)$ decreases instantaneously - jumps - changing
by $\epsilon^-<0$. If, after a positive jump, the value
of $x_i(t)$ exceeds the right boundary (i.e., $x_i(t)>1$), then a ``buy" is
considered to have occurred and the number of buys for that time interval is
increased by one.  Similarly, after $x_i(t)$ crosses the left boundary (i.e.,
$x_i(t)<-1$), the number of sells is increased by one. In either case,
$x_i(t)$ is set back to the neutral state (i.e., $x_i(t)=0$ with a small random offset to prevent
discontinuities) after the
number of buys or sells is updated. In this way, each individual agent's
decision affects the population's collective behavior.
This discrete jump process, combined with the previously described exponential decay, form the evolutionary rule for each agent's state $x_i(t)$,
which can be summarized as the following stochastic differential equation
\begin{equation}
                   dx_i = -\gamma x_i(t)dt + dI_i(t).
\label{xi:eqn}
\end{equation}

\subsection{The continuum model} \label{sec:cont}

It is possible to derive a concise approximate description of the dynamics of a large
assembly of agents by keeping track of only the density $\rho(x,t)$ of agents at each
state $x$, rather than the individual states of every agent in the
population. The resulting approximate evolution equation for $\rho(x,t)$, averaged over a large number
of replicas of the population, is given by \cite{Omurtag:01}
\begin{equation}
\frac{\partial \rho}{\partial t} =- \gamma \frac{\partial{(x\rho)}}{\partial x}
+\sum_{k=+,-}\nu^k(\rho(x-\epsilon^k,t)-\rho(x,t))+(R^++R^-)\delta(x).
\label{pop:eqn}
\end{equation}
For small $\epsilon^\pm$ by expanding $\rho(x-\epsilon^\pm,t)$ about
$\rho(x,t)$ in terms of $\epsilon^\pm$ and truncating terms higher than
second order in $\epsilon^\pm$, one obtains a Fokker-Planck-type approximation,
\begin{equation}
	\frac{\partial\rho}{\partial t}= \frac{\partial}{\partial x}(\mu\rho)
                                      +\frac{1}{2}\sigma^2\frac{\partial^2\rho}{\partial x^2}
                                      +(R^+ + R^-)\delta(x)\,,
\label{ori:eqn}
\end{equation}
where
\begin{equation}
	\mu (x,t)=\gamma x-(\nu^+\varepsilon^+ + \nu^-\varepsilon^-) \,,
\end{equation}
and
\begin{equation}
	 \sigma^2(t)=\nu^+(\varepsilon^+)^2+\nu^-(\varepsilon^-)^2 \,.
\end{equation}
Because agents leave the domain at $x=\pm 1$ and are restored at the origin, we have $\rho(x=\pm1,t)=0$, and a Dirac $\delta$ birth term at the origin. The buy and sell rates are defined as the outgoing fluxes at the boundaries,
$R^\pm=\mp\frac{1}{2}\sigma^2\partial\rho / \partial x\mid_{x=\pm1}$. The Fokker-Planck-type
equation~(\ref{ori:eqn}), with its non-local reinjection term, is the approximate continuum model
in our study.

\section{The finite volume scheme} \label{FVS:subsec}

We are now ready to describe a finite volume discretization of \e{ori:eqn}, which will serve as the first ``inner'' fine-scale simulator
on which our exploratory patch dynamics scheme (described in section 4) will be based. We divide
the one dimensional spatial domain into equal-sized grid cells and keep
track of our approximations to the average densities in these cells. At
each time step we update the average cell densities using approximations
of the fluxes through the edges of the cells. As shown in
\fg{fig:Scheme}, we denote the $i^{th}$ grid cell by
    \begin{equation}
	   \it{C}_i = (x_{i-1},x_i) \,;
    \end{equation}
\begin{figure}
\centering
{
    \includegraphics[scale=0.45]{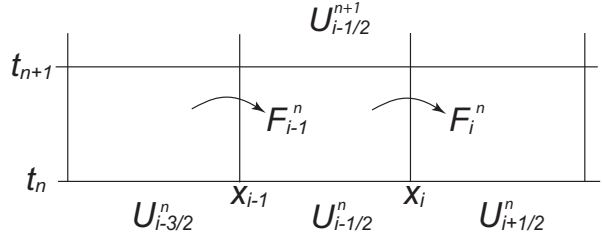}
}
\caption{Illustration of a finite volume scheme for updating the cell average $U_{i-1/2}^n$ by fluxes at the cell edges. }
\label{fig:Scheme}
\end{figure}
$U_{i-1/2}^n$ denotes the approximate average cell density in the $i$th cell at time $t_n$,
    \begin{equation}
	   U_{i-1/2}^n \approx \frac{1}{\Delta x}\int_{x_{i-1}}^{x_i}\rho(x,t_n)dx \equiv\frac{1}{\Delta
         x}\int_{\it{C}_i}\rho(x,t_n)dx\,,
    \label{Q:eqn}
    \end{equation}
where $\Delta x = x_{i}-x_{i-1}$ denotes the size of the cell.
Except for the central grid cell, which will be treated separately because of the particle reinjection, \eq{ori:eqn} can be rewritten as
\begin{equation}
	\frac{\partial\rho}{\partial t}= \frac{\partial}{\partial x}(\mu\rho
+\frac{1}{2}\sigma^2\frac{\partial\rho}{\partial x})\,.
    \label{rearrange:eqn}
\end{equation}
Define the function $f(x,t)$ as the flux term
    \begin{equation}
	   f(x,t)=- (\mu(x,t)\rho(x,t)+ \frac{1}{2}\sigma^2(t)\frac{\partial\rho(x,t)}{\partial x} ) \,,
    \label{f:eqn}
    \end{equation}
so that the integral form (in space) of \eq{rearrange:eqn} gives
    \begin{equation}
	   \frac{d}{dt}\int_{\it{C_i}} \rho(x,t)dx = f(x_{i-1},t) - f(x_i,t)  \,.
    \label{inSpace:eqn}
    \end{equation}
Integrating \eq{inSpace:eqn} in time from $t_n$ to $t_{n+1}$ and
dividing by $\Delta x$ gives
     \begin{equation}\nonumber
	   \frac{1}{\Delta x} \int_{\it{C}_i} \rho(x,t_{n+1})dx = \frac{1}{\Delta x} \int_{\it{C}_i} \rho(x,t_{n})dx - \frac{1}{\Delta x} \left[\int_{t_n}^{t_{n+1}}f(x_i,t)dt - \int_{t_n}^{t_{n+1}}f(x_{i-1},t)dt\right]  \,.
    \label{complete:eqn}
    \end{equation}
Based on \eq{Q:eqn}, \eq{complete:eqn} can be rewritten as
    \begin{equation}
	   U_{i-1/2}^{n+1} = U_{i-1/2}^n + \frac{\Delta t}{\Delta x}(F_{i-1}^n - F_i^n) \,,
    \label{overall:eqn}
    \end{equation}
    where $U_{i-1/2}^n$ is the average density defined in eq.(\ref{Q:eqn}),
    and $F_{i-1}^n$ and $F_i^n$  denote the
    average fluxes (i.e. fluxes averaged over $\Delta t$) across cell edges, that is,
    \eqn{flux:eqn1}
	   F_{i-1}^n = \frac{1}{\Delta t}\int_{t_n}^{t_{n+1}}f(x_{i-1},t)dt~~\rm{and}~~F_i^n = \frac{1}{\Delta t}\int_{t_n}^{t_{n+1}}f(x_i,t)dt.
    \enn
Based on (\ref{f:eqn}), the average fluxes over the time step are approximated by
     \begin{eqnarray}
    \nonumber  F_{i-1}^n & \approx & \it{F}(U_{i-3/2}^n,U_{i-1/2}^n) \\
                     &\approx & -( \mu_{i-1}\frac {U_{i-3/2}^n + U_{i-1/2}^n}{2}
                     + \frac{1}{2}\sigma^2 \frac {U_{i-1/2}^n - U_{i-3/2}^n}{\Delta x} ) \label{flux:eqn3}, \\
     \nonumber  F_{i}^n & \approx & \it{F}(U_{i-1/2}^n,U_{i+1/2}^n) \\
                     &\approx & -( \mu_{i}\frac {U_{i-1/2}^n + U_{i+1/2}^n}{2}
                     + \frac{1}{2}\sigma^2 \frac {U_{i+1/2}^n - U_{i-1/2}^n}{\Delta x} )                                 \,.
    \label{flux:eqn2}
    \end{eqnarray}
For a detailed discussion about the construction of diffusion fluxes for
finite volume schemes see \cite{DiffusionFlux}. The fluxes
at the outer two boundaries are expressed as (recall that $x_0=-1$ and $x_N=1$)
    \eqn{BC1:eqn}
	   F_0^n = -R^-  \approx  -\sigma^2 \frac {U_{1/2}^n}{\Delta x}~~\rm{and}~~F_N^n = R^+  \approx  \sigma^2 \frac {U_{N-1/2}^n}{\Delta x}.
    \enn
We use an odd number of grid cells $N$ for this scheme. Since the central
grid cell centered at $x_{N/2}$ contains the source point at $x=0$, two additional terms need to
be added to \eq{overall:eqn} for this cell,
    \begin{equation}
	   U_{N/2}^{n+1} = U_{N/2}^n + \frac{\Delta t}{\Delta x}(F_{(N-1)/2}^n - F_{(N+1)/2}^n + (R^+ + R^-)) \,.
    \label{central:eqn}
    \end{equation}
Clearly, the above scheme is conservative by construction. 	In section 5, this
	scheme is used as the ``inner" solver to mimic
	the agent-based simulator. Moreover, the associated patch dynamics scheme,
introduced in the next section, will be inspired by this scheme, and will also be conservative at the agent level.

\section{The patch dynamics scheme}\label{patchD}
\begin{figure}[th]
\centering
{
    \includegraphics[scale=0.45]{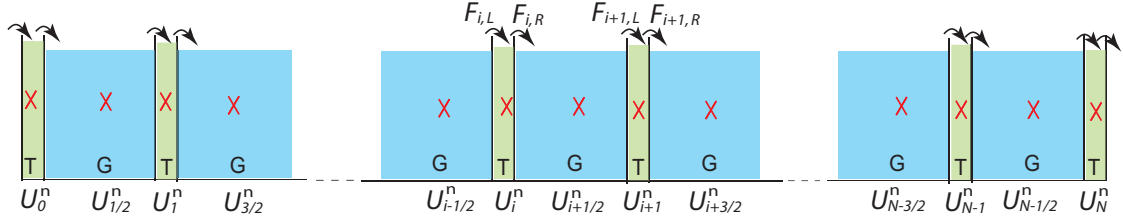}
}
\caption{Illustration of the gap-tooth scheme.}
\label{fig:Patch1}
\end{figure}
\begin{figure}[ht]
\centering
{
    \includegraphics[scale=0.45]{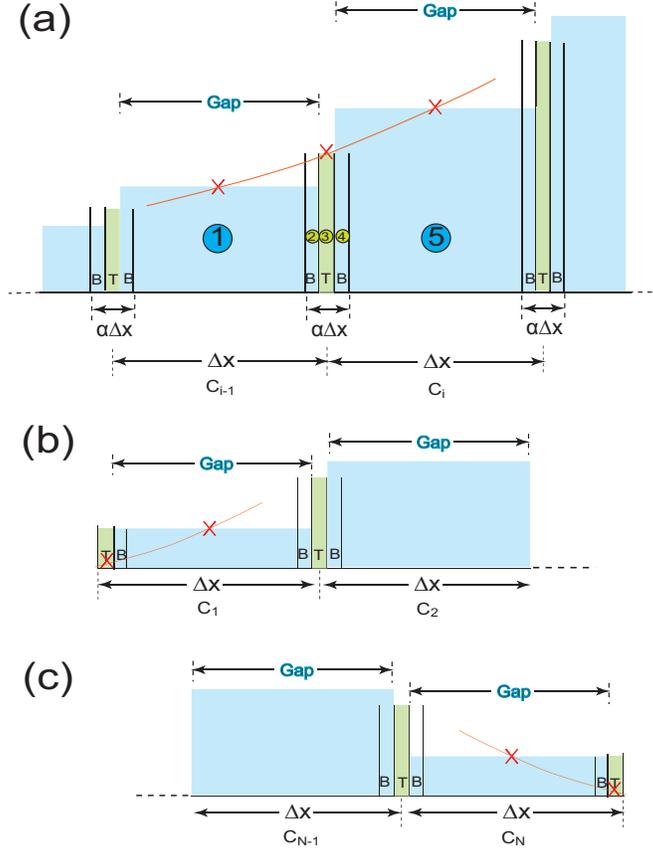}
}
\caption{Illustration of the lifting procedure. (a) lifting procedure
in the bulk of the domain. A local quadratic polynomial is constructed
such that the averaged polynomial values inside the tooth (region 3), and
inside the gaps (region 1 and 5) equal the average densities in these regions.
(b) lifting procedure in the left boundary. (c) lifting procedure in the right
boundary. One point value (zero on the the outer boundary) and two average values are used in constructing the
local polynomials for the simulation units at the two boundaries.}
\label{fig:Lifting_2}
\end{figure}

 As mentioned previously, the patch dynamics scheme combines gap-tooth and projective integration schemes. The gap-tooth scheme is
 illustrated in \fg{fig:Patch1}. We have divided the one dimensional
 space into two kinds of unequal-sized grid cells: the narrower
 cells, which we call teeth, and the wider cells, which we call
 gaps. We intend to solve the microscopic model only within the teeth (plus some buffer regions to handle the teeth boundary conditions). However, we keep track of the average densities in {\em both} the teeth and the gaps to obtain a conservative scheme. As in the finite volume scheme, we want to update the average cell
 densities based on the fluxes at the edges of the cells,
   \begin{eqnarray}
	   U_{i}^{n+\delta} &=& U_{i}^n + \frac{\delta t}{L_T}(F_{i,L}^n - F_{i,R}^n)   \,,    \label{Patch1:eqn}   \\
       U_{i+1/2}^{n+\delta} &=& U_{i+1/2}^{n} +  \frac{\delta t}{L_G}(F_{i,R}^n - F_{i+1,L}^n)  \,, \label{Patch2:eqn}
    \end{eqnarray}
where $L_T$ and $L_G$ denote the tooth and gap size respectively.
However, in contrast to the finite volume scheme, the fluxes
$F_{i,L}^n$ and $F_{i,R}^n$ are not computed based on a
known macroscopic equation (which is usually not available for
agent-based models). Instead, our goal is to compute these fluxes
based on the ``microscopic simulations" inside a small fraction of
the spatial domain -  the simulation unit (``$\alpha \Delta x$"
shown in Fig.~\ref{fig:Lifting_2}). These separated simulation
units are the only locations where microscopic simulations take place.
The numerical simulation of the microscopic model in each simulation
unit provides information on the evolution of the global state
 at the spatial location of the simulation unit, as if we were
 running the microscopic simulations in the entire spatial domain.
 Therefore, it is crucial to choose appropriate boundary conditions, so that the solutions inside
the teeth evolve as if the simulations were performed over the
entire domain. For some cases it is possible to implement
macroscopically-inspired constraints on the microscopic model
as boundary conditions \cite{Giovanni2005}. However, this is
not always the case. To overcome this difficulty, as discussed in \cite{Giovanni2006}, we use a larger box - the
simulation unit - to run the microscopic simulations, but still
compute the flux used to update the macroscopic state at the teeth boundaries. Each simulation unit consists of one tooth
and its edge buffers; the simulation units on the outer
boundaries of the spatial domain are treated slightly differently, as discussed later.
The purpose of the additional computational domains, the buffers,
is to ``protect'' the teeth simulations from boundary artifacts. This can be
accomplished over short enough time intervals, provided the
buffers are large enough.
As \fg{fig:Patch1} shows, we divide the one-dimensional
spatial domain into an odd number, $N$, of grid cells $\it{C_{i}}, i = 1,2,...,N$. The size of each cell is $\Delta x = 2/N$. We then put teeth (narrow bins) at the edges of these cells and put gaps
(wide bins) between the teeth. There are $\it{N+1}$ teeth and
$\it{N}$ gaps. The average densities for the teeth at time $t_n$
are denoted as $U_i^n, i=0,2,...N$, and the densities  for the
gaps are denoted as $U_{i-1/2}^n, i=1,2,...N$. We want to update
these macroscopic properties at time $t_{n+1}$ based on the
microscopic simulations inside the simulation units during the
time step from $t_n$ to $t_n+\delta t$.

The patch dynamics algorithm to proceed from $U_i^n$ to $U_i^{n+1}
(i = 0, \frac{1}{2}, 1, \cdots, N)$ is given below:
\\
\\
1. {\it Lifting}. At time $t_n$ create initial conditions
$u^i(x,t_n), i=0,...N$ inside each simulation unit for the
microscopic simulator, consistent with the spatial profile of
the macroscopic properties - the average densities inside the
teeth and gaps $U_i^n (i = 0, 1/2, 1, ...N)$. \\
\\
2. {\it Simulation}. Based on the microscopic initial condition
constructed in step 1 compute $F_{i,L}$ and $F_{i,R}$ by running
the microscopic simulator inside the simulation units from time
$t_n$ to $t_n + \delta t$.
\\
\\
3. {\it Restriction}. Obtain the spatially averaged densities
inside the teeth and gaps $U_i^{n+\delta} (i = 0, 1/2, 1, ...N)$  based on
Eq.~(\ref{Patch1:eqn})-(\ref{Patch2:eqn}).
\\
\\
4. {\it Projective step}. Estimate the macroscopic time derivative
at time $t_n$, e.g., as \\
 \begin{equation}
	   \mathcal{F}(U_i^n;\delta t) = \frac{U_i^{n+\delta}-U_i^n}{\delta t}     \,,
 \label{TimeDer:eqn}
 \end{equation}
and use it in a time integration method of choice, e.g., forward Euler,
\begin{equation}
	   U_i^{n+1} = U_i^n + \Delta t \mathcal{F}(U_i^n;\delta t) \,.
 \label{FEuler:eqn}
 \end{equation}

\subsection{Lifting Procedure} \label{sec:lifting}
Except for the two outer teeth that receive special treatment, we create the microscopic initial condition inside the $\it{i^{th}}$ simulation
unit by constructing a local quadratic polynomial approximation, $u^i(x,t)$,
to the density that matches the average density in a tooth and its two neighboring
gaps as follows.  Let $H$ be the width of the tooth (so $\Delta x - H$ is the
width of the gap) and require
\begin{eqnarray}
\frac{1}{H}\int_{-H/2}^{H/2} u^i(x,t_n)dx &=& U_i^n  \label{Eqn:Local1}, \\
\frac{1}{\Delta x - H} \int_{-\Delta x + H/2}^{-H/2}u^i(x,t_n)dx &=& U_{i-1/2}^n \label{Eqn:Local2}, \\
\frac{1}{\Delta x - H} \int_{H/2}^{\Delta x -H/2}u^i(x,t_n)dx &=& U_{i+1/2}^n \label{Eqn:Local3}.
\end{eqnarray}
As shown in Fig.~\ref{fig:Lifting_2}a, the microscopic initial conditions $u^i(x,t_n)$
$(i=1,2,...N-1)$ are constructed such that the averaged function values inside
the tooth (region 3), and inside the neighboring gaps (regions 1 and 5) equal the average
densities in these regions.

At the two boundaries, as Fig.~\ref{fig:Lifting_2}b and c show, there
are gaps only on \emph{one} side of the boundary teeth so the microscopic initial
conditions are created in a slightly different manner.
 At the left boundary,
a local quadratic polynomial $u^0(x,t_n)$ is constructed such that
      \begin{eqnarray}
      u^0(-1,t_n) &=& 0,  \\
      \frac{1}{H}\int_{-H/2}^{H/2} u^0(x,t_n)dx &=& U_0^n,  \\
      \frac{1}{\Delta x - 3H/2} \int_{H/2}^{\Delta x -H}u^0(x,t_n)dx &=& U_{1/2}^n.
      \end{eqnarray}
Similarly, at the right boundary a local quadratic polynomial $u^N(x,t_n)$
is constructed such that
  \begin{eqnarray}
      u^N(1,t_n) &=& 0,  \\
      \frac{1}{H}\int_{-H/2}^{H/2} u^N(x,t_n)dx &=& U_N^n,  \\
      \frac{1}{\Delta x - 3H/2} \int_{-\Delta x +H}^{-H/2}u^N(x,t_n)dx &=& U_{N-1/2}^n.
  \end{eqnarray}
Note that the \emph{left edge} of the leftmost tooth (and similarly the
right edge of the rightmost tooth) is put at the two boundaries $x=\pm1$,
instead of their respective centers. For this reason, the outermost
gaps are of a slightly different size ($\Delta x - \frac{3}{2}H$) compared to
the inner ones (which are of equal size $\Delta x - H$).

\subsection{Simulation step} \label{sec:simulation}

As mentioned previously, for our initial experiment and its analysis a fine-grid finite volume scheme is first used as
the``inner", microscopic simulator. We emphasize again that the purpose
of first using a fine discretization of \e{ori:eqn} as our ``inner" microscopic solver is to facilitate a
numerical study of the errors involved - something we
cannot do when dealing directly with the agent-based simulator.

\begin{figure}
\centering
{
    \includegraphics[scale=0.73]{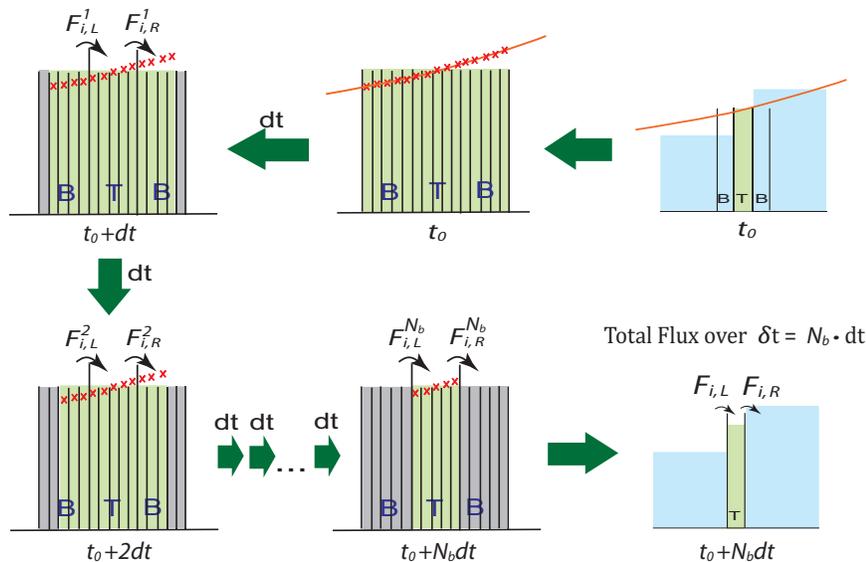}
}
\caption{Illustration of the simulation step for the simulation units in the bulk of the spatial domain.}
\label{fig:Simulation}
\end{figure}

\begin{figure}
\centering
{
    \includegraphics[scale=0.52]{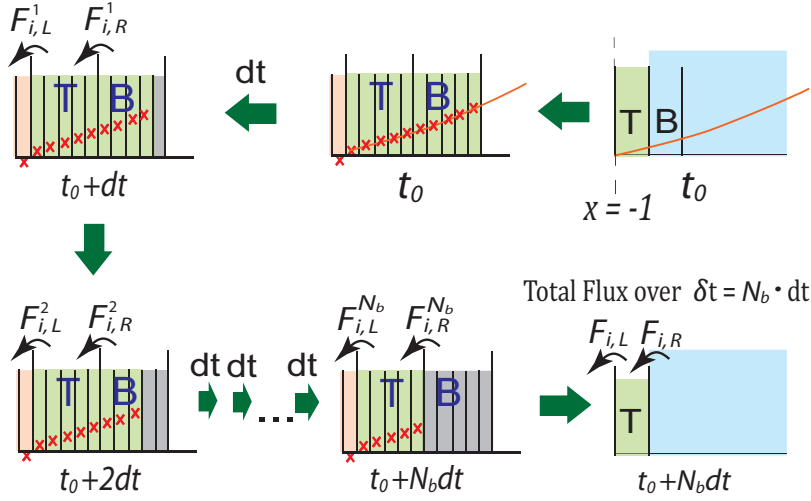}
}
\caption{Illustration of the simulation step for the simulation unit at the left boundary.}
\label{fig:SimulationBC}
\end{figure}

We divide each simulation unit into $N_t+2N_b$
fine bins, where $N_t$ denotes the number of fine bins inside the tooth
while $N_b$ denotes the number of fine bins inside each buffer. All
these fine bins are of equal size $\delta x$. Based on the previously constructed
local quadratic polynomial $u^i(x,t_n)$ we initialize the averaged densities
in the fine bins as follows:
\begin{equation}
	   \frac{1}{h}\int_{x_j-h/2}^{x_j+h/2} u^i(x,t_n) dx = u_j^i \,,\\
\label{eqn:finebin}
\end{equation}
where $x_j$ denotes the center of each fine bin and $u_j^i$ denotes the
average densities inside the $j^{th}$ fine bin of the $i^{th}$ simulation
unit (as the red cross marks in Fig.~\ref{fig:Simulation} and
\ref{fig:SimulationBC} denote). Fig.~\ref{fig:Simulation} illustrates the
simulation step in the bulk of the spatial domain. Starting from the
previously constructed local quadratic polynomial, we first initialized
the average densities inside the fine bins using \eq{eqn:finebin}.
To evolve $u_j^i(t_0)$ for one microscopic step $dt$, as Eq. (\ref{overall:eqn})
and (\ref{flux:eqn3}-\ref{flux:eqn2}) show, it requires $u_{j-1}^i(t_0)$,
$u_{j}^i(t_0)$ and $u_{j+1}^i(t_0)$. Because of that, we are not able to evolve
the solutions for the two bins \emph{at the boundaries} (the two bins denoted in grey
in the upper-left part of Fig.~\ref{fig:Simulation}). Therefore, after the
first microscopic time step $dt$, these two bins are discarded. For the same
reason, after another microscopic time step dt, two more (now outer) fine bins need to be
discarded, and so on until all the bins in the buffer regions are
discarded. There are $N_b$ fine bins in each buffer, so we can run the
microscopic simulator for $N_b$ microscopic time steps. During the $k$-th  step ($k=1,2,...N_b$), we save the fluxes
computed at the left (resp.~right) edge of the $i$-th tooth, $F_{i,L}^k$ (resp.~ $F_{i,R}^k$).

The simulations in the leftmost and rightmost simulation units are performed
slightly differently. As Fig.~\ref{fig:SimulationBC} shows, a
buffer is only used on one side of the tooth (the side to the interior of the spatial
domain). On the other side, the outer boundary of the spatial domain, the boundary condition $\rho(x=\pm1)=0$ needs to be enforced.  This is achieved using a ghost bin (colored pink in Fig.~\ref{fig:SimulationBC}) with density
equal in  magnitude but opposite in sign to
the density of the leftmost (resp. rightmost) inner bin. (See Fig.~\ref{fig:SimulationBC} for an illustration of the left boundary;
the right boundary is treated in a similar way.)

\subsection{Restriction procedure}

At the end of the $N_b^{th}$ step, we update the \emph{total} fluxes at the
edges of each tooth for the coarse time interval $\delta t = N_b\cdot dt$
as follows
\begin{eqnarray}
F_{i,L} &=& \frac{1}{N_b}\sum_{k=1}^{N_b}F_{i,L}^k, \\
F_{i,R} &=& \frac{1}{N_b}\sum_{k=1}^{N_b}F_{i,R}^k.
\end{eqnarray}
Then, based on Eqn.~(\ref{Patch1:eqn} - \ref{Patch2:eqn}), the coarse
variables - the average densities inside the teeth and
gaps $U_i^{n+\delta} (i = 0, 1/2, 1, ...N)$ - are updated.

\subsection{Projective step} \label{sec:project}

\eq{TimeDer:eqn} is now used to estimate the macroscopic
time derivative and \eq{FEuler:eqn} is used to project the
macroscopic solutions forward to obtain $U_i^{n+1}$.

\section{Results and discussion}

\subsection{Comparison of the patch-dynamics scheme with fine and coarse finite volume schemes for the continuum equation}
We compare the performance of the patch-dynamics scheme with two other
schemes: one is the fine grid finite volume scheme i.e., our
microscopic simulator {\em applied to the entire region}; the other is a coarse grid finite volume scheme, which mimics the usually unavailable macroscopic solver.
All three schemes are first applied to the continuum approximation
of the agent-based model.
\begin{figure}[h]
\centering{
    \includegraphics[scale=0.52]{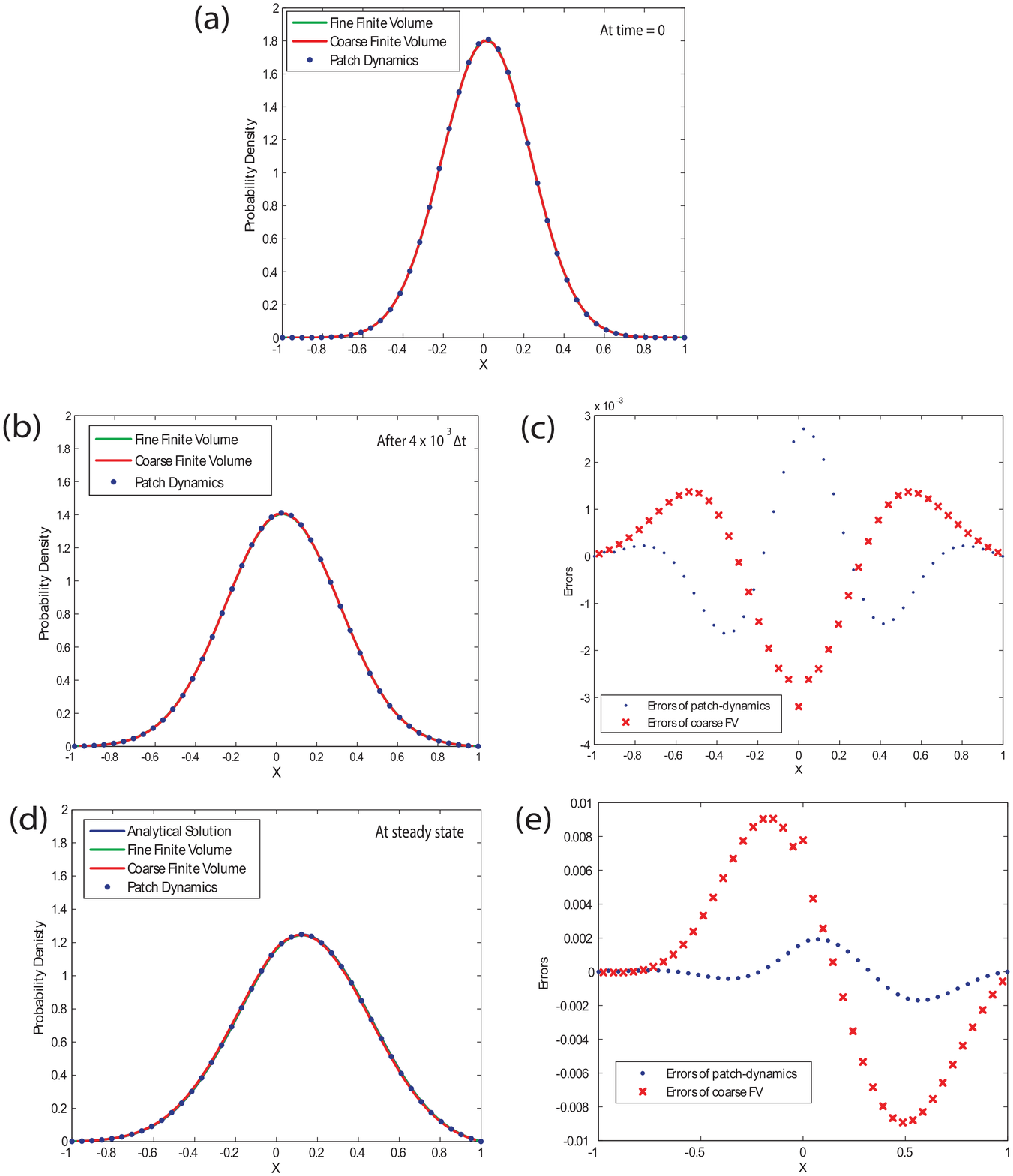}
    }
\caption{Comparison of the solutions obtained using the three different schemes. (a) The initial conditions
for the three schemes, which are constructed from a common Gaussian-like function
$f(x) = a \cdot exp\left(-\left(\left(x-b\right)/c\right)^2\right)$ with $a = 1.811,
b = 0.01545, c = 0.3115$.
(b) Solutions of the three schemes at time $4\times 10^3\Delta t$.
(c) Comparison of the errors of the patch-dynamics scheme and the coarse finite volume scheme at
time $4\times 10^3\Delta t$, using the fine finite volume scheme as the reference.
(d) Steady state solutions.
(e) Comparison of the errors of the patch-dynamics scheme and the coarse finite volume scheme at steady state.
The blue curve (in (a), (b), and (e)) denotes the solution of the fine finite volume scheme;
the red curve (in (a), (b), and (e)) denotes the solution of the coarse finite volume scheme; the blue
dots (in (a), (b), and (e)) denote the solutions of the teeth for the patch dynamics scheme;
the blue dots (in (c) and (e)) denote the errors (in teeth) of the patch dynamics scheme;
The red cross marks in (in (c) and (e)) denote the errors of the coarse finite volume scheme.
Parameter values for the three schemes are as follows: $g = 40,
\epsilon^+ = 0.075, \epsilon^- = -0.072, \gamma = 1, \nu_{ex}^+= \nu_{ex}^-=20$.
Parameters specific for the patch dynamics scheme: $\alpha = 0.1, N = 41,
N_b = N_t = 10, dt = 4h^2, \delta t = 9dt, M = 90, \Delta t = (M+1)\delta t$.
Parameters specific to the coarse finite volume scheme: $N = 41, \Delta t=
(M+1)\delta t$. Parameters specific to the fine finite volume scheme:
$N_{fine} = 1271$.}
\label{fig:Evolution}
\end{figure}
Fig. \ref{fig:Evolution} shows several snapshots of the solutions
of these three schemes along the time-path towards their steady states.
Starting from the same Gaussian-like initial conditions (Fig. \ref{fig:Evolution}a)
we evolve the solutions of the three systems for a long time
(Fig. \ref{fig:Evolution}b shows one snapshot of the transient solutions)
until they have reached their steady states (Fig. \ref{fig:Evolution}d).
As we can see from the figure, the solution of the patch-dynamics scheme
agrees well with the other two schemes along the entire process. The steady
state solutions of the three schemes also agree well with the analytical
steady state solution of the problem (Fig. \ref{fig:Evolution}d).
In Fig. \ref{fig:Evolution}c and \ref{fig:Evolution}e we plot the differences between
the fine finite-volume scheme and the other two methods, patch dynamics and coarse finite volume scheme. The solution to the fine finite-volume scheme is an approximation to the true solution, so these differences can be viewed as a first approximation to the errors in the other two methods.
Also, since the patch dynamics uses the fine finite-volume scheme as its microscopic solver, the differences between it and patch dynamics are an indication of the error introduced by the patch dynamics scheme.
As we can see from Fig. \ref{fig:Evolution}c and \ref{fig:Evolution}e, during the transient phase
the two schemes give comparable errors, while as we further evolve the system to steady state the patch dynamics scheme gives smaller errors. This is probably due to the fact that our patch dynamics scheme utilizes more microscopic level information.

We now turn to the computational savings of patch dynamics.
For the patch-dynamics scheme, we divide the spatial domain into $N=41$ big cells, giving
42 teeth and 41 gaps. We run ``microscopic simulations" in only $\alpha=10\%$ of
the spatial domain (since the width of each simulation unit is equal to
only $10\%$ of each big cell). Inside each simulation unit we put $10$
equal-sized fine bins inside each tooth and each buffer to run the microscopic simulations
(mimicked here by the fine grid finite volume scheme).
We use $dt = 4\times h^2$ as the time step
for the microscopic simulator ($\delta x$ is the width of the ``fine" bin).
After running the microscopic
simulator for $n_b = 10$ microscopic steps (dt), or equivalently one coarse
time step $\delta t = n_b \times dt$, we project macroscopic variables 90 coarse steps
forward in time.
This means that the microscopic simulations are performed in only $1.11\%$ of the
temporal domain.
For the fine finite volume scheme, we use $N_{fine} = 1271$
fine bins, and use the same time step size dt as the one used in the
microscopic simulator for the patch dynamics scheme.
For the coarse
finite volume scheme we use the same number of big cells ($N=41$) as the one
used for the patch-dynamics scheme.
The time step size of this scheme
is set to be equal to the effective time step size $\Delta t$ of the
patch dynamics scheme (recall that $\Delta t = (M+1)\delta t$, with $M = 90$).
Overall, in this illustrative example the patch dynamics scheme runs the
microscopic simulations in only $10\%$ of the spatial domain and $1.11\%$ of the temporal domain.

\subsection{Error Analysis}
The patch-dynamics scheme uses a microscopic integrator that is not specified (in this paper we have used an agent-based model as well as a fine-scale finite volume method to illustrate the method and to permit better analysis).  The purpose of the gap-tooth part of the scheme is to estimate time derivatives of the macroscopic tooth and gap densities (more precisely, the chords of their time-dependent solutions over a small interval $\Delta t$), and these are input to an unspecified outer integrator or other macroscopic numerical analysis procedure (in this paper we have used projective forward Euler integration, but any scheme could be used).

The errors thus arise from three sources: the microscopic integrator itself, the patch-dynamics scheme, and the outer integrator (or whatever other numerical procedure we apply to the chord estimates).  The final error is a non-linear combination of the errors from all three sources that depends on details of the problem being solved and the details of each method, making it very difficult to isolate each part of the process.
Hence it seems advisable to consider each step from a {\it backward error} perspective.
In the backward error view, the microscopic integrator integrates a {\it perturbed} equation exactly.  The size of the perturbation is a function of the microscopic integrator and the equation being integrated, but for our purposes we can ignore it and assume that the microscopic integrator is {\em exact} (it is for the perturbed equation).  Then we can examine the size of the errors introduced by the gap-tooth scheme.

The equations for estimating the teeth and gap chord slopes (for other than the center and end teeth) are:
\begin{equation}
	 \frac{U_{i}^{n+1} - U_{i}^n}{\Delta t} = \frac{F_{i,L}^n - F_{i,R}^n}{H} \label{eqn:Utooth}
\end{equation}
and
\begin{equation}
     \frac{U_{i-1/2}^{n+1} - U_{i-1/2}^{n}}{\Delta t} = \frac{F_{i-1,R}^n - F_{i,L}^n}{\Delta x - H} \,. \label{eqn:Ugap}
\end{equation}
Note that if we compute the fluxes $\{F_{i,L}^n\}$ and $\{F_{i,R}^n\}$ exactly, these equations are exact.  If, for example, we ran the microscopic integrator over all space and time, these equations would give the exact values of the average teeth and gap densities over time.  However, we introduced the teeth to avoid running the microscopic simulation over all of space, so that after a finite time the buffer regions fail to ``protect" the teeth from the lack of simulation information from the bulk of the gaps.  At that point we have to stop the microscopic integration, and lift from the macroscopic information using the process described in \s{patchD} to get a new microscopic description.  Under the assumption that the microscopic integrator is correct, the error introduced by the patch dynamics scheme is simply {\em a change of initial values} at the start of the $\Delta t$ integration block.

In the next subsection we will examine this error analytically and then report on some computational estimates of the errors in our finance model in the following subsection.

\subsubsection{Lifting Error}

Suppose that the computed microscopic solution at $t_n$ is $v(x)$.  For each tooth the lifting process creates a quadratic polynomial approximation, $u(x)$, such that the average densities of $u$ and $v$ agree in the tooth and its neighboring gaps.  The microscopic integration then continues starting from the $u$ values rather than from $v$. The effect is to add
\eqn{pdef}
r(x) = u(x) - v(x)
\enn
to the initial values for the next integration interval.

We analyze each tooth separately.  Let us set the spatial origin, $x = 0$, to be the tooth center and assume that the Taylor series for $v(x)$ is
\eqn{taylorv}
v(x) = v_0 + v_1x + v_2x^2/2 + v_3x^3/6 + \cdots ,
\enn
i.e., $v_q$ is the $q$-th spatial derivative of $v(x)$.
It will be convenient to define $h = H/2$ and $D = \Delta x - h$ (so the left and right gaps are $[-D, -h]$ and $[h,D]$.  Because of the average density condition in the tooth and its surrounding gaps we have
\begin{eqnarray}
\int_{-D}^{-h} r(x)dx & = & 0, \nonumber \\
\int_{-h}^{h} r(x)dx & = & 0, \label{intcond} \\
\int_{h}^{D} r(x)dx & = & 0, \nonumber
\end{eqnarray}
where $r(x)$ is given in \e{pdef}.
Let us define $u(x)$ as
\eqn{udef}
u(x) = r_0 + v_0 + (r_1+v_1)x + (r_2 + v_2)x^2/2.
\enn
Hence
\eqN
r(x) = u(x) - v(x) = r_0 + r_1x + r_2x^2/2 - v_3x^3/6 - v_4x^4/24 + \cdots
\enN
Solving \e{intcond} for $r_i$ in terms of $v_j$ we get
\begin{eqnarray}
r_0 & = & -D^2h^2v_4/120 + {\rm O}(D^6 + h^6)\label{esol0}, \\
r_1 & = & (D^2 + h^2)v_3/24 + {\rm O}(D^4 + h^4)\label{esol1}, \\
r_2 & = & (D^2 + h^2)v_4/20  + {\rm O}(D^4 + h^4)\label{esol2}.
\end{eqnarray}
Thus, over the region comprising the tooth and its neighboring gaps the largest error that can be introduced is
\eqn{maxerr}
E = \max_{x \in [-D, D]} |[(D^2+h^2)x -4x^3]v_3/24| + {\rm O}(h^4 + D^4)
\enn
which is ${\rm O}(h^3 + D^3)$.

How does the error in Eq. \ref{maxerr} affect the solution of the integration problem?
 Unfortunately it depends on the problem, but if the problem is such that modest perturbations do not cause serious difficulties (for example, if the system after semi-discretization in space satisfies a Lipschitz condition) the final error is equal to the sum of the errors, each multiplied by bounded amplifications.  Thus, if the number of steps in the outer integrator is ${\rm o}((\max(D,h)^{-3})$, the global error will be bounded.  If the average time step size is ${\rm O}((\max(D,H))^2)$ the global error is ${\rm O}(\max(D,H))$ over any finite interval, and if the average time step size is ${\rm O}(\max(D,H))$ the global error is ${\rm O}((\max(D,H))^2)$.  If the solution tends to a stable stationary state, then the final error may reflect only the most recent lifting errors and if many of the components are stiff, this may be true over much of the integration interval.

\subsubsection{Numerical Results}

\begin{figure}
\centering
{
    \includegraphics[scale=0.5]{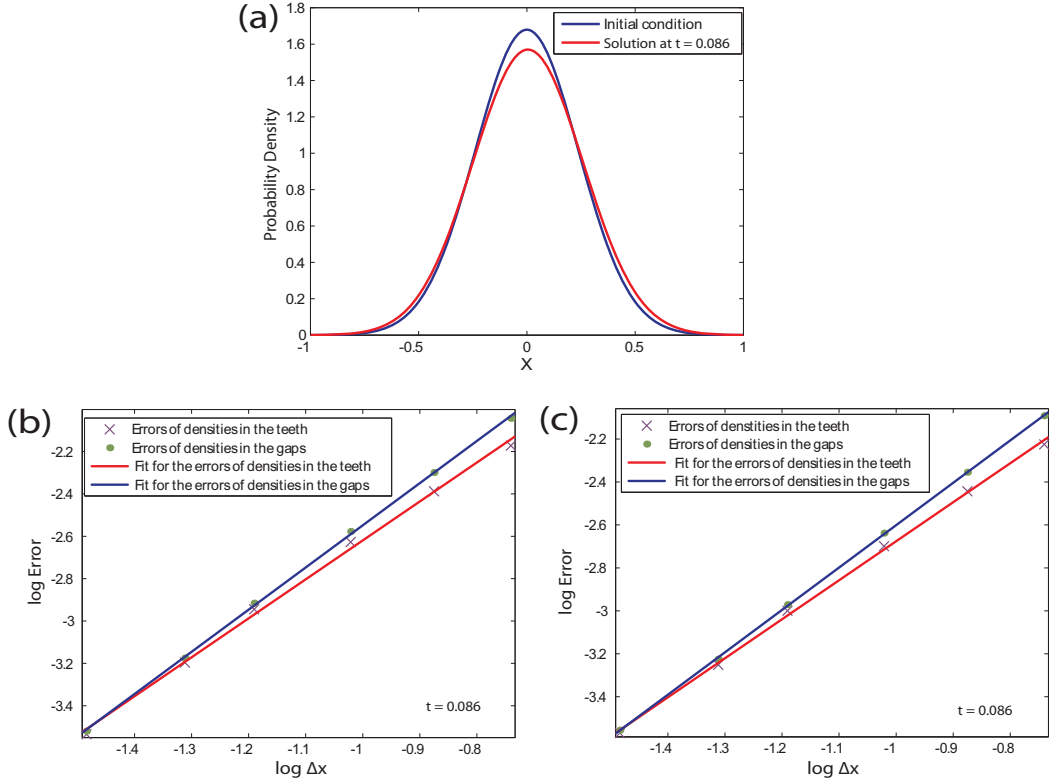}
}
\caption{(a) The blue curve corresponds to the Gaussian-like initial
condition for both the patch-dynamics scheme and the reference scheme.
The red curve corresponds to the solution of the reference scheme at
$t = 0.086$. This was obtained by running the reference scheme for $10000$
time steps. The patch-dynamics scheme was running for the same amount of
time. The differences of the average densities inside the teeth as well as the gaps
between these two schemes were taken  at this time point to generate the
log-log plots. The model parameters used in these simulations are the same
as the ones listed in Fig. \ref{fig:Evolution}. (b) Log(Error) v.s. log ($\Delta x$)
plot for $\emph{Limit 1}$. (b) Log(Error) v.s. log ($\Delta x$) plot for 
$\emph{Limit 2}$. In (b) and (c), the dots denote the log(Error) v.s.
log($\Delta x$) data points for the teeth, while the cross marks denote the
corresponding data points for the gaps. The blue and red curves denote the
linear fit for the data sets of the teeth and gaps respectively. The R-square
values for the fits are all very close to one (0.998 and above). The slopes
of the fitted lines (order of consistency estimates) together with their
respective $95\%$ confidence intervals are presented in Table \ref{tab:Order}.}
\label{fig:Error}
\end{figure}

To numerically validate that the order of consistency for the patch
dynamics scheme is two, we compare solutions of the patch-dynamics
scheme (using a fine-grid finite volume scheme as the ``inner" microscopic simulator)
with the solutions of a reference scheme, which is also a fine-grid finite volume
scheme (but over the entire physical domain).
In this reference scheme, the spatial domain is divided into $963$
fine bins with bin width $\delta x = 2/963$. The microscopic time step is
set to be $dt = 2\delta x^2$. As Fig. \ref{fig:Error}a shows, starting from the
same Gaussian-like initial condition, we evolve both the patch-dynamics
scheme and the reference scheme for some time ($t=0.086$, or equivalently
$10000$ time steps in the reference scheme), and then compute the differences
between the solutions of these two schemes in order to construct the
subsequent log-log plot for the consistency order analysis.

For the patch-dynamics scheme, two distinct limits have been numerically simulated.
The first limit corresponds
to cases where the tooth size becomes very small. In practice, the tooth
size in this limit is set to be equal to the fine bin size in the
reference scheme ($\delta x = 2/963$). The size of the buffers is set to
be the same as the one of the teeth. By fixing the size of the teeth
and buffers, we studied cases with different sizes of big cells and
therefore different sizes of gaps. Recall that, as Fig. \ref{fig:Lifting_2}
shows, each big cell is of size $\Delta x$, which is equal to the sum of
the sizes of a pair of tooth and gap. In the second limit, instead of
fixing the teeth size to have a very small value, we keep the ratio
between the sizes of the simulation unit and the big cell
fixed, while varying both of them simultaneously. The second limit is
closely related to the envisioned practical simulation cases, because in practice
we need to use teeth with some finite size, in order to capture
meaningful macroscopic properties.

For each limit we compute the differences between the average
densities inside the teeth and the gaps for the patch dynamics
scheme on the one hand, and the average densities inside the corresponding regions
of the reference scheme on the other. We define the errors for the teeth and gaps as the
L2 norm of these differences:
\begin{eqnarray}
\left \| U_t-U_{ref} \right \| &=& \sqrt{\frac{1}{N+1} \sum_{i=0}^{N}(U_t^i-U_{ref}^i)^2}  \label{eqn:norm1}  \\
\left \| U_g-U_{ref} \right \| &=& \sqrt{\frac{1}{N} \sum_{j=1/2}^{N-1/2}(U_t^j-U_{ref}^j)^2} \label{eqn:norm2}
\end{eqnarray}
where, $U_t^i$ and $U_g^j$ denote the average densities in the
teeth and gaps for the patch-dynamics scheme, while $U_{ref}^i$
and $U_{ref}^j$ denote the average densities inside the corresponding
regions of the reference scheme. The detailed results of the different
cases in these two limits are shown in Table \ref{tab:NumLimits}.

To estimate the order of consistency, we generated the log(error)
v.s. log($\Delta x$) plot and performed linear fitting.
The fitted plots for the two limits are shown in Fig. \ref{fig:Error} (b)
and (c). The slope of the fitted lines together with their respective $95\%$
confidence interval estimates are shown in Table \ref{tab:Order}.
The numerical consistency order estimates are all close to $2$;
the slight deviation is most likely due to the discontinuity at the origin
and the boundary conditions.
%

\begin{table}[t]
{\renewcommand{\arraystretch}{1.5}
 \renewcommand{\tabcolsep}{0.2cm}
  \begin{tabular}{ | c | c | c | c | c | c | c |}
    \hline
      \multicolumn{2}{|c|}{} & \multicolumn{2}{|c|}{Limit 1} & \multicolumn{2}{|c|}{Limit 2} \\ \hline
    N & $\Delta x$  & $\left \| U_t-U_{ref} \right \|$ & $\left \| U_g-U_{ref} \right \|$  & $\left \| U_t-U_{rEf} \right \|$ & $\left \| U_g-U_{rEf} \right \|$  \\ \hline
      $ 11 $ & $1.82E-1$ & $6.76E-03$ & $9.12E-03$ & $6.03E-03$ & $8.13E-03$
      \\ \hline
      $ 15 $ & $1.33E-1$ & $4.07E-03$ & $5.01E-03$ & $3.63E-03$ & $4.47E-03$
      \\ \hline
      $ 21 $ & $9.52E-2$ & $2.29E-03$ & $2.63E-03$ & $2.04E-03$ & $2.29E-03$
      \\ \hline
      $ 31 $ & $6.45E-2$ & $1.15E-03$ & $1.20E-03$ & $1.00E-03$ & $1.07E-03$
      \\ \hline
      $ 41 $ & $4.88E-2$ & $6.61E-04$ & $6.92E-04$ & $5.75E-04$ & $6.17E-04$
      \\ \hline
      $ 61 $ & $3.28E-2$ & $2.88E-04$ & $3.02E-04$ & $2.69E-04$ & $2.75E-04$
      \\ \hline
  \end{tabular}}
  \caption{Numerical results of the consistency order analysis for
  the two limits. N is the number of big cells in the patch-dynamics scheme.
  $\Delta x$ denotes the size for the big cells. The error terms are defined
  in Eq. (\ref{eqn:norm1}-\ref{eqn:norm2}). }
  \label{tab:NumLimits}
\end{table}

\begin{table}[t]
{\renewcommand{\arraystretch}{1.5}
 \renewcommand{\tabcolsep}{0.2cm}
  \begin{tabular}{ | c | c | c | }
    \hline
       & Limit 1 & Limit 2 \\ \hline
    $U_t$ & 1.839 (1.715, 1.964) & 1.818 (1.732, 1.903) \\ \hline
    $U_g$ & 1.991 (1.939, 2.043) & 1.974 (1.948, 2.001) \\ \hline
  \end{tabular}}
  \caption{Consistency order estimates for the two limits. $95\%$ confidence intervals for the order estimates are shown in the brackets.}
  \label{tab:Order}
\end{table}

\subsection{Patch dynamics scheme for the agent-based model} \label{agent:sec}
Finally, we turn to the agent-based model
 (discussed in section \ref{sec:agent}) for which the patch dynamics scheme
 has been designed. We divide the
 spatial domain into $N=21$ big cells, and set the ratio of the
 size of the simulation unit to the big cell ($\alpha$) to 0.2.
 This means we are running the microscopic simulations in
 $20\%$ of the spatial domain. Our macroscopic observables (our coarse variables)
 are the average agent densities inside the teeth and gaps. Inside each
 simulation unit the tooth and the two buffers are of equal size,
 and each of them has been subdivided into 10 fine bins. Based on the
 lifting procedure discussed in section \ref{sec:lifting}, microscopic
 agent densities are assigned to these fine bins. We then convert these
 density values into numbers of agents based on the total amount of
 agents available in the system; in this case the total number of agents we use is $N_{agents} = 3\times 10^{6}$.
 Once the number of agents for each fine bin is assigned we
 start the agent-based simulations. As illustrated in section
 \ref{sec:simulation}, the coarse time step size $\delta t$ - the
 time we can run the agent-based simulator before we stop for reconstruction -
 is closely related to the size of the buffers. This is because the
 buffers are gradually ``contaminated'' due to boundary artifacts;
 we need to stop a little before these artifacts get propagated to the tooth.
 Our coarse time step size $\delta t$ is chosen
 to be $2\times 10^{-3}$. Within this time interval, the chance for any
 agent to travel a distance of the size of a buffer or more is very low.
 In this way agents outside each simulation unit will not affect
 the solution of the tooth (although they do affect the solution of the
 buffers), and the tooth is therefore protected. Along the simulation process,
 the fluxes of agents across the edges of the teeth are tracked.
 In section \ref{sec:simulation}, the fluxes are updated based on
 equations of the artificial microscopic simulator (a PDE). In the
 agent-based case we track the fluxes by counting the number of
 agents going across those edges.

Due to the stochastic nature of the agent-based modeling, we also
need to reduce the noise for the macroscopic variables. The
agent-based simulations are therefore repeated inside each simulation
unit for $N_{realizations} = 20$ copies before the restriction procedure
is applied. At the end of each coarse time step, we perform the
restriction procedure to compute the averaged fluxes at the edges
of the teeth. Based on these fluxes we update the numbers of agents
inside the teeth as well as the neighboring gaps, and compute the corresponding
densities based on the numbers of agents. Following the same procedure
discussed in section \ref{sec:project}, based on the macroscopic
solutions (the densities) at the start and the end of each coarse time
step, we estimate local time derivatives of the macroscopic variables
and project them forward in time. What allows us to do this is the
smoothness of these coarse variables in the time domain. In this case,
we run the agent-based simulations over the duration of one coarse time step, and
then ``jump" nine coarse time steps ahead. In other words, we are running the agent-based
simulations over $10\%$ of the temporal domain.

\begin{figure}
\centering{\includegraphics[scale=0.6]{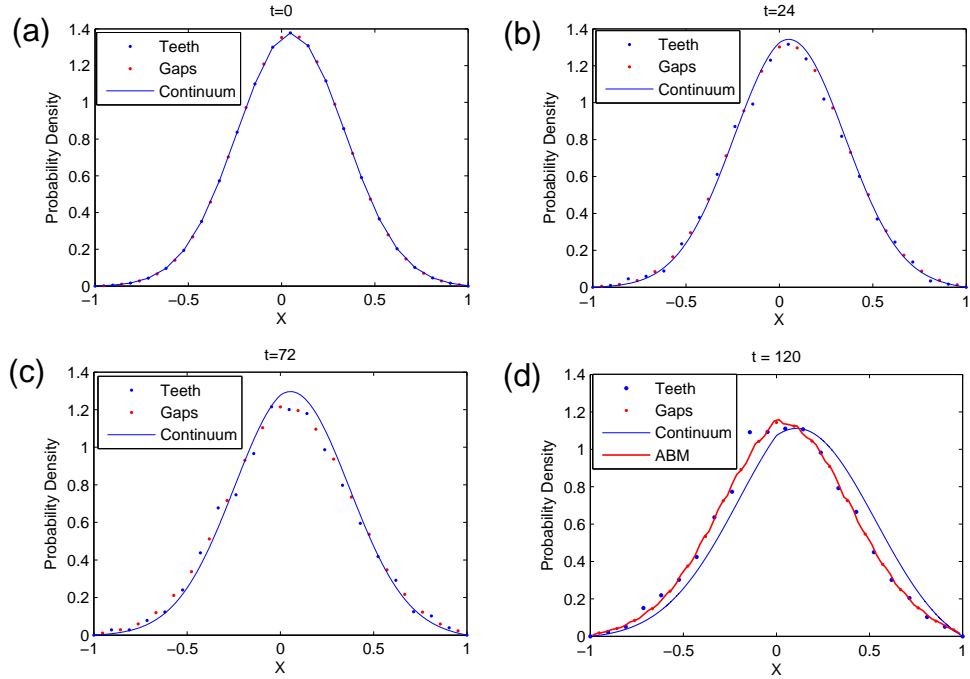}}
\caption{Snapshots along a patch dynamics agent-based trajectory
on its way to an eventual steady state. The average agent densities inside
the teeth are denoted by red dots, while the ones inside the gaps
are denoted by blue dots. To compare the performance, solutions of a fine-grid ($N=629$)
finite volume scheme for the continuum model are plotted in (a-c) in blue curves.
In (d) the analytical steady state of the continuum model is plotted as a blue curve while
the approximated steady state of the full agent-based simulations is plotted
as a red curve. Parameter values for the patch dynamics scheme
 are as follows: $N=21$, $\alpha = 0.2$,  $g = 1$, $\epsilon^+ = 2.9\times 10^{-3}$,
$\epsilon^- = -2.89 \times 10^{-3}$, $\gamma = 1\times 10^{-3}$,
$N_{agents} = 3\times 10^{6}$, $N_{realizations} = 20$,
$\delta t = 2\times 10^{-3}$, $M = 9$.}
\label{fig:Agent}
\end{figure}

Fig. \ref{fig:Agent} shows four snapshots along a patch dynamics
agent-based trajectory on its way to the eventual steady state. The solutions for the
agent-based patch-dynamics simulations are plotted in blue (teeth) and red (gaps)
respectively. To compare the performance, the solutions of a fine-grid finite volume scheme for the
continuum model are also plotted (in blue curve)as a reference. The analytical steady
state for the continuum model is shown in Fig. \ref{fig:Agent}(d).
Overall, the solutions of the agent-based patch-dynamics simulations match well
(although not perfectly) with the solutions of the continuum model. The deviations
are most likely due to the approximations made to derive the continuum model from the agent-based one.
In Fig. \ref{fig:Agent}(d) we can see that the steady state of the patch-dynamics agent
based simulations matches the stationary state of the full scale agent-based simulations better than
the analytical steady state of the continuum model matches it.

We reiterate that the patch-dynamics agent-based simulations
require the ``expensive" agent-based computations to be performed in only $10\%$ of the temporal domain and
$20\%$ of the spatial domain compared to the full scale agent-based
simulations.

\section{Conclusions}
We described a patch dynamics scheme for conservative agent-based problems.
This scheme approximates an unavailable effective equation over
macroscopic time and length scales, when only a microscopic agent-based
simulator is given. It only uses appropriately initialized simulations
of the agent-based model over small subsets (patches) of the spatiotemporal
domain, thus significantly reducing the computational cost. Because this
scheme mimics a finite volume scheme for the underlying effective equation,
it is conservative by construction. Since it is often not possible to
impose macroscopically inspired boundary conditions on a microscopic
agent-based simulation, we have used buffer regions around the patches,
which temporarily shield the internal region of the patches from boundary
artifacts. We have explored numerical characteristics of the scheme based on the
continuum approximation (a Fokker-Planck-type evolution equation) of the agent-based model, and demonstrated its
effectiveness for agent-based computations involving
a financial market agent-based model. In this paper the agent-based model
was mainly used for illustration purposes, to demonstrate the effectiveness of the approach.
Several factors that can affect the performance of the patch dynamics scheme still need to be
explored in more detail. One of the most interesting ones to study should be the effect
of the noise (in the teeth and the gaps) on the estimation of relevant quantities,
and through them to the overal performance of this scheme.

{\bf Acknowledgements.} This work was partially supported by the US Department of Energy
and the National Science Foundation.

\newpage


\end{document}